\documentclass[11pt,leqno,fleqn]{article}
\usepackage{amsmath,amsfonts,hyperref}
\newcommand{\dps}{\displaystyle}
\newcommand{\dy}[2]{%
\refstepcounter{equation}%
\LABEL{#1}%
\begin{list}{}{
\topsep 5mm
\leftmargin 18mm
\rightmargin 0cm
\itemsep 0mm
\listparindent 0mm
\parsep 0mm
\itemsep 0mm
\labelsep 0mm
\labelwidth 18mm
}%
\item[\rm (\theequation)\hfill]
#2
\end{list}%
}
\newcommand{\dyz}[1]{%
\refstepcounter{equation}%
\begin{list}{}{
\topsep 5mm
\leftmargin 18mm
\rightmargin 0cm
\itemsep 0mm
\listparindent 0mm
\parsep 0mm
\itemsep 0mm
\labelsep 0mm
\labelwidth 18mm
}%
\item[\rm (\theequation)\hfill]
#1
\end{list}%
}
\newcommand{\dyyz}[1]{\dyz{\raggedright$\dps#1$}}
\newcommand{\dyy}[2]{\dy{#1}{\raggedright$\dps#2$}}
\newcommand{\de}[2]{\dy{#1}{\raggedright$\displaystyle #2 $}}
\newcommand{\dez}[1]{\dyz{\raggedright$\displaystyle #1 $}}
\newcounter{sectie}
\newcounter{subsectie}
\newcommand{\sectz}[1]{\refstepcounter{sectie}\setcounter{subsectie}{0}
\section*{\boldmath \thesectie. #1}%
}
\newcommand{\sect}[2]{\refstepcounter{sectie}\setcounter{subsectie}{0}
\section*{\boldmath \thesectie. #2}%
\label{#1}}
\newcommand{\oN}{{\mathbb{N}}}
\newcounter{stelling}
\newcommand{\thm}[2]{\refstepcounter{stelling}\vspace{4mm}\noindent{\bf Theorem \thestelling.}\label{#1}{\it #2}}

\newcommand{\corz}[2]{\vspace{4mm}\noindent{\bf Corollary \thestelling #1.}{\it #2}}
\newcommand{\oR}{{\mathbb{R}}}
\newcommand{\FF}{{\cal F}}
\newcommand{\LABEL}[1]{\label{#1}}
\newcommand{\rf}[1]{{\rm (\ref{#1})}}
\newcommand{\subsect}[2]{\refstepcounter{subsectie}
\subsection*{\boldmath \thesectie.\thesubsectie. #2}%
\label{#1}}
\newcounter{bewering}
\newcommand{\prop}[2]{\refstepcounter{bewering}\vspace{4mm}\noindent{\bf Proposition \thebewering.}\label{#1}{\it #2}}
\newcommand{\pf}{\vspace{3mm}\noindent{\bf Proof.}\ }
\newcommand{\kfrac}[2]{\mbox{$\frac{#1}{#2}$}}
\newcommand{\bx}{\hspace*{\fill} \hbox{\hskip 1pt \vrule width 4pt height 8pt depth 1.5pt \hskip 1pt}

\addvspace{4mm}}
\newcommand{\bxx}{\hspace*{\fill} \hbox{\hskip 1pt \vrule width 4pt height 8pt depth 1.5pt \hskip 1pt}}
\begin{document}

\begin{center}
{\LARGE\bf Compact orbit spaces in Hilbert spaces and limits of edge-colouring models

}
\end{center}
\vspace{1mm}
\begin{center}
{\large
\hspace{10mm}
Guus Regts\footnote{ University of Amsterdam.
The research leading to these results has received funding from the European Research Council
under the European Union's Seventh Framework Programme (FP7/2007-2013) / ERC grant agreement
n$\mbox{}^{\circ}$ 339109.}
and
Alexander Schrijver$\mbox{}^1$

}

\end{center}

\noindent
{\small{\bf Abstract.}
Let $G$ be a group of orthogonal transformations of a real Hilbert space $H$.
Let $R$ and $W$ be bounded $G$-stable subsets of $H$.
Let $\|.\|_R$ be the seminorm on $H$ defined by
$\|x\|_R:=\sup_{r\in R}|\langle r,x\rangle|$ for $x\in H$.
We show that if $W$ is weakly compact and
the orbit space $R^k/G$ is compact for each $k\in\oN$,
then the orbit space $W/G$ is compact when $W$ is equiped with
the norm topology induced by $\|.\|_R$.

As a consequence we derive the existence of limits of edge-colouring models
which answers a question posed by Lov\'asz.
It forms the edge-colouring counterpart of the graph limits of
Lov\'asz and Szegedy, which can be seen as limits of vertex-colouring models.

In the terminology of de la Harpe and Jones, vertex- and edge-colouring models
are called `spin models' and `vertex models' respectively.

}

\sectz{Introduction}

In a fundamental paper, Lov\'asz and Szegedy [10]
showed that the collection of simple graphs fits in a natural
way in a compact metric space $W$ that conveys several phenomena
of extremal and probabilistic graph theory and of
statistical mechanics.
In particular, a limit behaviour of graphs can be derived.

The elements of $W$ are called {\em graphons}, as generalization of
graphs, but they can also be considered to generalize the vertex-colouring
models, or `spin models' in the sense of de la Harpe and Jones [7].
In this context, the partition functions of
spin models form a compact space.
In the present paper, we investigate to what extent
edge-colouring models, or `vertex models' in the terminology
of [7], behave similarly.
Indeed, the edge-colouring models form a dense subset in a compact space,
and thus we obtain limits of edge-colouring models.
This solves a problem posed by Lov\'asz [8].

To obtain these results, we prove a general theorem on compact orbit spaces in
Hilbert space, that applies both to vertex- and to edge-colouring models.
This compactness theorem uses and extends theorems of
Lov\'asz and Szegedy [11]
on Szemer\'edi-like regularity in Hilbert spaces.

For background on graph limits we also refer to the recent book of Lov\'asz [9].
Partition functions of edge-colouring models with a finite
number of states were characterized by
Szegedy [14]
and
Draisma, Gijswijt, Lov\'asz, Regts, and Schrijver [2].

\sect{9jl12f}{Formulation of results}

In this section we describe our results, giving proofs in subsections
\ref{9jl12f}.\ref{9jl12b},
\ref{9jl12f}.\ref{9jl12c}, and
\ref{9jl12f}.\ref{9jl12e}.
Throughout, for any Hilbert space $H$, $B(H)$ denotes the closed
unit ball in $H$.
Moreover, measures are Lebesgue measures.

\bigskip
\noindent
{\bf Compact orbit spaces in Hilbert spaces.}
Let $H$ be a real Hilbert space and let $R$ be a bounded subset of $H$.
Define a seminorm $\|.\|_R$ and a pseudometric\footnote{A
{\em seminorm} is a norm except that nonzero elements may have norm 0.
A {\em pseudometric}
is a metric except that distinct points may have distance 0.
One can turn a pseudometric space into a metric space by identifying points
at distance 0, but for our purposes it is notationally easier and
sufficient to maintain the original space.
Notions like convergence pass easily over to pseudometric spaces, but
limits need not be unique.}
$d_R$ on $H$ by, for $x,y\in H$:
\dyyz{
\|x\|_R:=\sup_{r\in R}|\langle r,x\rangle|
\text{~~~and~~~}
d_R(x,y):=\|x-y\|_R.
}
In this paper, we use the topology induced by this
pseudometric only if we explicitly mention it,
otherwise we use the topology
induced by the usual Hilbert norm $\|.\|$.

A subset $W$ of $H$ is called {\em weakly compact}
if it is compact in the weak topology on $H$.
Note that for any $W\subseteq H$ (cf.\ [6]):
\dyz{
$W$ closed, bounded, and convex $\Rightarrow$
$W$ weakly compact $\Rightarrow$
$W$ bounded.
}

Let $G$ be a group acting on a topological space $X$.
Then the {\em orbit space} $X/G$ is the quotient space of $X$ taking
the orbits of $G$ as classes.
A subset $Y$ of $X$ is called
{\em $G$-stable} if $g\cdot y\in Y$ for each $g\in G$ and $y\in Y$.

Our first main theorem (which we prove in Section \ref{9jl12f}.\ref{9jl12b}) is:

\thm{13me12a}{
Let $H$ be a Hilbert space and
let $G$ be a group of orthogonal transformations of $H$.
Let $W$ and $R$ be $G$-stable subsets of $H$, with $W$ weakly compact
and $R^k/G$ compact for each $k$.
Then $(W,d_R)/G$ is compact.
}

\bigskip
\noindent
{\bf Application to graph limits and vertex-colouring models.}
As a consequence of Theorem \ref{13me12a} we now first derive the compactness of
the graphon space, which was proved by Lov\'asz and Szegedy [11].
Let $H=L^2([0,1]^2)$, the set of all square integrable functions $[0,1]^2\to\oR$.
Let $R$ be the collection of functions $\chi^A\times\chi^B$, where $A$ and
$B$ are measurable subsets of $[0,1]$ and where $\chi^A$ and $\chi^B$ denote
the incidence functions of $A$ and $B$.
Let $S_{[0,1]}$ be the group of measure space automorphisms of $[0,1]$.
The group $S_{[0,1]}$ acts naturally on $H$.
Moreover, $R^k/S_{[0,1]}$ is compact for each $k$.
(This can be derived from the fact that for each measurable $A\subseteq[0,1]$
there is a $\pi\in S_{[0,1]}$ such that $\pi(A)$ is an interval up to a
set of measure 0 (cf.\ [12]).)

Let $W$ be the set of $[0,1]$-valued functions $w\in H$ satisfying $w(x,y)=w(y,x)$ for all $x,y\in[0,1]$.
Then $W$ is a closed bounded convex $S_{[0,1]}$-stable subset of $H$.
So by Theorem \ref{13me12a}, $(W,d_R)/S_{[0,1]}$ is compact.
The elements of $W$ are called {\em graphons},
where two elements $w,w'$ of $W$ are assumed to
be the same graphon if $w'=g\cdot w$ for
some $g\in S_{[0,1]}$.
Therefore one may say that the graphon space is
compact with respect to $d_R$.

In the context of
de la Harpe and Jones [7], graphons can be considered as `spin models'
(with infinitely many states).
For any $w\in W$, the {\em partition function} $\tau(w)$ of $w$
is given by, for any graph $F$:
\dez{
\tau(w)(F):=\int_{[0,1]^{V(F)}}\prod_{uv\in E(F)}w(x(u),x(v))dx.
}
Let $\FF$ denote the collection of simple graphs.
Lov\'asz and Szegedy [10] showed that $\tau:(W,d_R)\to\oR^{\FF}$
is continuous (here the restriction to simple graphs is necessary).
Since $(W,d_R)/S_{[0,1]}$ is compact and since $\tau$ is
$S_{[0,1]}$-invariant, the collection of functions $f:\FF\to\oR$ that are
partition functions of graphons is compact.
Hence each sequence $\tau(w_1),\tau(w_2),\ldots$ of partition functions
of graphons $w_1,w_2,\ldots$,
such that $\tau(w_i)(F)$ converges for each $F$,
converges to the partition function $\tau(w)$ of some graphon $w$.

Since simple graphs can be considered as elements of $W$ (by considering their
adjacency matrix as element of $W$), this also gives a limit behaviour of
simple graphs.

\bigskip
\noindent
{\bf Application to edge-colouring models.}
We next show how Theorem \ref{13me12a} applies to
the edge-colouring model (also called vertex model).
For this, it will be convenient to use a different, but universal
model of Hilbert space.
Let $C$ be a finite or infinite set, and
consider the Hilbert space $H:=\ell^2(C)$, the set of all functions
$f:C\to\oR$ with $\sum_{c\in C}f(c)^2<\infty$,
having norm $\|f\|:=(\sum_{c\in C}f(c)^2)^{1/2}$.
(Each Hilbert space is isomorphic to $\ell^2(C)$ for some $C$.)
Following de la Harpe and Jones [7], any element of $\ell^2(C)$ is called
an {\em edge-colouring model}, with {\em state set} (or {\em colour set}) $C$.

Define for each $k=0,1,\ldots$:
\dez{
H_k:=\ell^2(C^k).
}
The tensor power $\ell^2(C)^{\otimes k}$ embeds naturally in $\ell^2(C^k)$,
and the orthogonal group $O(H)$ of $H$ acts naturally on $H_k$.
Define moreover
\dez{
R_k:=\{r_1\otimes\cdots\otimes r_k\mid r_1,\ldots,r_k\in B(H)\}
\subseteq H_k.
}

Again, let $\FF$ be the collection of simple graphs.
As usual, $H_k^{S_k}$ denotes the set of elements of $H_k$ that are invariant under the natural action of $S_k$ on $H_k$.
Define the function
\dyy{19jl12a}{
\pi:\prod_{k=0}^{\infty}H^{S_k}_k\to\oR^{\FF}
\text{~~~by~~~}
\pi(h)(F):=
\sum_{\phi:E(F)\to C}\prod_{v\in V(F)}h_{\deg(v)}(\phi(\delta(v)))
}
for $h=(h_k)_{k\in\oN}\in\prod_{k\in\oN}H_k^{S_k}$
and $F\in\FF$.
Here $\deg(v)$ denotes the degree of $v$.
Moreover, $\delta(v)$ is the set of edges incident
with $v$, in some arbitrary order, say
$e_1,\ldots,e_k$, with $k:=\deg(v)$.
Then $\phi(\delta(v)):=
(\phi(e_1),\ldots,\phi(e_k))$ belongs to $C^k$.
(So $\phi(\delta(v))$ may be viewed as the set
of colours `seen' from $v$.)
For \rf{19jl12a}, the order chosen is irrelevant, as $h_k$
is $S_k$-invariant.

The function $\pi(h):\FF\to\oR$ can be considered as the
{\em partition function} of the edge-colouring model $h$.
It is not difficult to see that $\pi$ is well-defined,
and continuous if we
take the usual Hilbert metric on each $H_k$, even if we replace
$\FF$ be the collection of all graphs without loops
(cf.\ \rf{18me12a}).
For simple graphs it remains continuous on $\prod_kB_k$
where $B_k:=B(H_k^{S_k})$
if we replace for each $k$ the Hilbert metric on $B_k$ by $d_{R_k}$:

\thm{17me12d}{
$\pi$ is continuous on $\dps\prod_{k\in\oN}(B_k,d_{R_k})$.
}

\bigskip
\noindent
This is proved in Section \ref{9jl12f}.\ref{9jl12c}, while in
Section \ref{9jl12f}.\ref{9jl12e} we derive from Theorem \ref{13me12a}:

\thm{17me12c}{
$\dps(\prod_{k=0}^{\infty}(B_k,d_{R_k}))/O(H)$ is compact.
}

\bigskip
\noindent
Now $\pi$ is $O(H)$-invariant.
This follows from the facts that $\ell^2(C^k)$
is the completion of $\ell^2(C)^{\otimes k}$
and that $O(H)$-invariance is direct if each
$h_k$ belongs to $\ell^2(C)^{\otimes k}$.
Hence Theorem \ref{17me12c} implies:

\corz{a}{
The image $\pi(\prod_kB_k)$ of $\pi$ is compact.
}

\bigskip
\noindent
This implies:

\corz{b}{
Let $h^1,h^2,\ldots\in\prod_kB_k$ be such that for each simple graph
$F$, $\pi(h^i)(F)$ converges as $i\to\infty$.
Then there exists $h\in\prod_kB_k$ such that for each simple graph
$F$,
$\lim_{i\to\infty}\pi(h^i)(F)=\pi(h)(F)$.
}

\bigskip
\noindent
As $\ell^2(C')$ embeds naturally in $\ell^2(C)$ if $C'\subseteq C$,
all edge-colouring models with any finite number of states embed in $\ell^2(C)$
if $C$ is infinite.
So the corollary also describes a limit behaviour of finite-state
edge-colouring models, albeit that the limit may have infinitely many states.

The corollary holds more generally for sequences in $\prod_k\lambda_kB_k$,
for any fixed sequence $\lambda_0,\lambda_1,\ldots\in\oR$.

We do not know if the quotient function
$\pi/\sim:(\prod_kB_k)/\sim\to \oR^{\FF}$ is one-to-one,
where $\sim$ is the equivalence relation on $\prod_kB_k$ with
$h\sim h'$ whenever $h'$ belongs to the closure of the $O(H)$-orbit
of $h$.
(For finite $C$ and $\FF$ replaced by the set of all graphs,
this was proved in [13].)
The analogous result for vertex-colouring models (i.e., graph limits) was
proved by Borgs, Chayes, Lov\'asz, S\'os, and Vesztergombi [1].

\subsect{9jl12b}{Proof of Theorem \ref{13me12a}}

In this section we give a proof of Theorem \ref{13me12a}.

\prop{28ja12e}{
Let $H$ be a Hilbert space and let $R,W\subseteq H$ with
$R$ bounded and $W$ weakly compact.
Then $(W,d_R)$ is complete.
}

\pf
Let $a_1,a_2,\ldots\in W$ be a Cauchy sequence with respect to $d_R$.
We must show that it has a limit in $W$, with respect to $d_R$.
As $W$ is weakly compact, it has the Bolzano-Weierstrass property;
that is, there is a point $a$ such that each weak neighbourhood of $a$
contains infinitely many terms of $a_1,a_2,\ldots$.
Then $a$ is a required limit, that is,
$\lim_{n\to\infty}d_R(a_n,a)=0$.
For choose $\varepsilon>0$.
As $a_1,a_2,\ldots$ is a Cauchy sequence with respect to $d_R$,
there is a $p$ such that $d_R(a_n,a_m)<\kfrac12\varepsilon$ for
all $n,m\geq p$.
For each $r\in R$,
as each weak neighbourhood of $a$ contains infinitely many terms of
$a_1,a_2,\ldots$, there is an $m\geq p$ with
$|\langle r,a_m-a\rangle| < \kfrac12\varepsilon$.
This implies for each $n\geq p$:
\dyyz{
|\langle r,a_n-a\rangle|\leq
|\langle r,a_n-a_m\rangle|+|\langle r, a_m-a\rangle|
<\varepsilon.
}
So $d_R(a_n,a)\leq\varepsilon$ if $n\geq p$.
\bx

Let $G$ be a group acting on a pseudometric space $(X,d)$ that leaves
$d$ invariant.
Define a pseudometric $d/G$ on $X$ by, for $x,y\in X$:
\dez{
(d/G)(x,y)=\inf_{g\in G}d(x,g\cdot y).
}
As $d$ is $G$-invariant, $(d/G)(x,y)$ is equal to the distance of the $G$-orbits
$G\cdot x$ and $G\cdot y$.
Any two points $x,y$ on the same $G$-orbit have $(d/G)(x,y)=0$.
If we identify points of $(X,d/G)$ that are on the
same orbit, the topological space obtained is homeomorphic to
the orbit space $(X,d)/G$ of the topological space $(X,d)$.
Note that being compact or not is invariant under such identifications.

\prop{2fe12a}{
If $(X,d)$ is a complete metric space, then $(X,d/G)$ is complete.
}

\pf
Let $a_1,a_2,\ldots\in X$ be Cauchy with respect to $d/G$.
Then it has a subsequence $b_1,b_2,\ldots$
such that $(d/G)(b_k,b_{k+1})<2^{-k}$ for all $k$.
Let $g_1=1\in G$.
If $g_k\in G$ has been found, let $g_{k+1}\in G$ be such that
$d(g_kb_k,g_{k+1}b_{k+1})<2^{-k}$.
Then $g_1b_1,g_2b_2,\ldots$ is Cauchy with respect to $d$.
Hence it has a limit, $b$ say.
Then $\lim_{k\to\infty}(d/G)(b_k,b)=0$, implying
$\lim_{n\to\infty}(d/G)(a_n,b)=0$.
\bx

Let $H$ be a Hilbert space and let $R\subseteq H$.
For any $k\geq 0$, define
\dez{
Q_k=\{\lambda_1r_1+\cdots+\lambda_k r_k\mid r_1,\ldots,r_k\in R,
\lambda_1,\ldots,\lambda_k\in[-1,+1]\}.
}

For any pseudometric $d$, let $B_d(Z,\varepsilon)$ denote the
set of points at distance at most $\varepsilon$ from $Z$.
The following is a form of `weak Szemer\'edi regularity'.
(cf.\ Lemma 4.1 of Lov\'asz and Szegedy [11],
extending results of
Frieze and Kannan [5] and
Fernandez de la Vega, Kannan, Karpinski, and Vempala [4] --- the method is essentially that of
[5]):

\prop{13ja12b}{
If $R\subseteq B(H)$, then for each $k\geq 1$:
\de{19jl12c}{
B(H)\subseteq B_{d_R}(Q_k,1/\sqrt{k}).
}
}

\pf
Choose $a\in B(H)$ and set $a_0:=a$.
If, for some $i\geq 0$,
$a_i$ has been found, and $d_R(a_i,0)>1/\sqrt{k}$,
choose $r\in R$ with $|\langle r,a_i\rangle|>1/\sqrt{k}$.
Define $a_{i+1}:=a_i-\langle r,a_i\rangle r$.
Then, with induction,
\dyy{7ok12c}{
\|a_{i+1}\|^2
=
\|a_i\|^2-\langle r,a_i\rangle^2(2-\|r\|^2)
\leq
\|a_i\|^2-\langle r,a_i\rangle^2
\leq
\|a_i\|^2-1/k
\leq
1-(i+1)/k.
}
Moreover,
as $\langle r,a_i\rangle\in[-1,+1]$, we know by induction that $a-a_i\in Q_i$.

By \rf{7ok12c}, as $\|a_{i+1}\|^2\geq 0$,
the process terminates for some $i\leq k$.
For this $i$ one has $d_R(a_i,0)\leq 1/\sqrt{k}$.
Hence, as $Q_i\subseteq Q_k$,
$d_R(a,Q_k)\leq d_R(a,Q_i)\leq d_R(a,a-a_i)=d_R(a_i,0)\leq 1/\sqrt{k}$.
\bx

\bigskip
\noindent
{\bf Proof of Theorem \ref{13me12a}:}
Observe that $R$ is bounded as $R/G$ is compact.
So we can assume that $R,W\subseteq B(H)$.

By Propositions \ref{28ja12e} and \ref{2fe12a}, $(W,d_R/G)$ is complete.
So it suffices to show that $(W,d_R/G)$ is
totally bounded; that is, for each
$\varepsilon>0$, $W$ can be covered by finitely
many $d_R/G$-balls of radius $\varepsilon$
(cf. [3]).

Set $k:=\lceil 4/\varepsilon^2\rceil$.
As $R^k/G$ is compact, $Q_k/G$ is compact
(since the function
$R^k\times[-1,1]^k\to Q_k$ 
mapping
$(r_1,\ldots,r_k,\lambda_1,\ldots,\lambda_k)$
to
$\lambda_1r_1+\cdots+\lambda_kr_k$ is
continuous, surjective, and $G$-equivariant).
Hence (as $d_R\leq d_2$) $(Q_k,d_R)/G$ is compact,
and so $(Q_k,d_R/G)$ is compact.
Therefore,
$Q_k\subseteq B_{d_R/G}(F,1/\sqrt{k})$ for some finite $F$.
Then with Proposition \ref{13ja12b},
\dez{
W\subseteq B(H)\subseteq B_{d_R}(Q_k,1/\sqrt{k})\subseteq B_{d_R/G}(Q_k,1/\sqrt{k})\subseteq B_{d_R/G}(F,2/\sqrt{k})\subseteq B_{d_R/G}(F,\varepsilon).
\bxx
}

\subsect{9jl12c}{Proof of Theorem \ref{17me12d}}

For any graph $F$, define a function
\de{10jl12a}{
\pi_F:\prod_{v\in V(F)}B_{\deg(v)}\to\oR
\text{~~~by~~~}
\pi_F(h):=\sum_{\phi:E(F)\to C}\prod_{v\in V(F)}h_v(\phi(\delta(v)))
}
for $h=(h_v)_{v\in V(F)}\in\prod_{v\in V(F)}B_{\deg(v)}$.
(The sum in \rf{10jl12a} converges, as follows from \rf{18me12a} below.)

\prop{14me12a}{
If $F$ is a simple graph, then $\pi_F$ is continuous on
$\dps\prod_{v\in V(F)}(B_{\deg(v)},d_{R_{\deg(v)}})$.
}

\pf
We first prove the following.
For any $k$, any $h\in H_k^{S_k}$, and any $c_1,\ldots,c_l\in C$ with $l\leq k$,
let $h(c_1,\ldots,c_l)$ be the element of $H_{k-l}^{S_{k-l}}$ with
$h(c_1,\ldots,c_l)(c_{l+1},\ldots,c_k)=h(c_1,\ldots,c_k)$ for
all $c_{l+1},\ldots,c_k\in C$.
Let $k_1,\ldots,k_n\in\oN$, let $h_i\in H_{k_i}^{S_{k_i}}$ for $i=1,\ldots,n$,
and let $F=([n],E)$ be a graph with $\deg(i)\leq k_i$ for
$i=1,\ldots,n$.
Then
\de{18me12a}{
\sum_{\phi:E\to C}\prod_{v\in [n]}\|h_v(\phi(\delta(v)))\|\leq
\prod_{v\in [n]}\|h_v\|.
}
We prove this by induction on $|E|$, the case $E=\emptyset$
being trivial.
Let $|E|\geq 1$, and choose an edge $ab\in E$.
Set $E':=E\setminus\{ab\}$
and $\delta'(v):=\delta(v)\setminus\{ab\}$
for each $v\in V(F)$.
Then
\dyyz{
\sum_{\phi:E\to C}\prod_{v\in [n]}\|h_v(\phi(\delta(v)))\|
=
\sum_{\phi:E'\to C}
\sum_{c\in C}\|h_a(\phi(\delta'(a)),c)\|\|h_b(\phi(\delta'(b)),c)\|
\prod_{v\in [n]\atop v\neq a,b}\|h_v(\phi(\delta(v)))\|
\leq
\sum_{\phi:E'\to C}
\|h_a(\phi(\delta'(a)))\|\|h_b(\phi(\delta'(b)))\|
\prod_{v\in [n]\atop v\neq a,b}\|h_v(\phi(\delta(v)))\|
\leq
\prod_{v\in [n]}\|h_v\|,
}
where the inequalities follow from Cauchy-Schwarz and induction,
respectively.
This proves \rf{18me12a}.

We next prove that for each $h=(h_v)_{v\in V(F)}\in\prod_{v\in V(F)}H_{\deg(v)}$
and each $u\in V(F)$:
\de{13me12b}{
\pi_F(h)\leq
\|h_u\|_{R_{\deg(u)}}
\prod_{v\in V(F)\atop v\neq u}\|h_v\|.
}
To see this,
let $N(u)$ be the set of neighbours of $u$,
$F':=F-u$, and $\delta'(v):=\delta(v)\setminus\delta(u)$ for $v\in V(F)\setminus\{u\}$.
Then
\dyyz{
\pi_F(h)
=
\sum_{\phi:E(F)\to C}\prod_{v\in V(F)}h_v(\phi(\delta_F(v)))
=
\sum_{\phi:E(F')\to C}
\langle \bigotimes_{v\in N(u)}h_v(\phi(\delta'(v))),h_u\rangle
\prod_{v\in V(F')\setminus N(u)}
h_v(\phi(\delta(v)))
\leq
\sum_{\phi:E(F')\to C}
\|h_u\|_{R_{\deg(u)}}
\prod_{v\in V(F')}
\|h_v(\phi(\delta'(v)))\|
\leq
\|h_u\|_{R_{\deg(u)}}
\prod_{v\in V(F')}\|h_v\|,
}
where the inequalities follow from the definition of $\|.\|_{R_{\deg(u)}}$ and
from \rf{18me12a} (applied to $F'$), respectively.
This proves \rf{13me12b}.

Now let $g,h\in\prod_{v\in V(F)}B_{\deg(v)}$.
We can assume that $V(F)=[n]$.
For $u=1,\ldots,n$, define $p^u\in\prod_{i\in[n]}B_{\deg(i)}$
by $p^u_i:=g_i$ if $i<u$, $p^u_u:=g_u-h_u$, and $p^u_i:=h_i$ if $i>u$.
Moreover, for $u=0.\ldots,n$, define $q^u\in\prod_{i\in[n]}B_{\deg(i)}$
by $q^u_i:=g_i$ if $i\leq u$ and $q^u_i=h_i$ if $i>u$.
So $q^n=g$ and $q^0=h$.
By the multilinearity of $\pi_F$ we have
$\pi_F(q^u)-\pi_F(q^{u-1})=\pi_F(p^u)$.
Hence by \rf{13me12b} we have the following, proving the theorem,
\dez{
\pi_F(g)-\pi_F(h)
=
\sum_{u=1}^n(\pi_F(q^u)-\pi_F(q^{u-1}))
=
\sum_{u=1}^n\pi_F(p^u)
\leq
\sum_{u=1}^n\|p^u_u\|_{R_{\deg(u)}}
=
\sum_{u=1}^n\|g_u-h_u\|_{R_{\deg(u)}}.
\bxx
}

\noindent
{\bf Proof of Theorem \ref{17me12d}:}
For each $F\in\FF$,
the function $\psi:\prod_{k\in\oN}B_k\to\prod_{v\in V(F)}B_{\deg(v)}$
mapping
$(h_k)_{k\in\oN}$ to $(h_{\deg(v)})_{v\in V(F)}$ is continuous.
As $\pi(.)(F)=\pi_F(\psi(.))$,
the theorem follows from Proposition \ref{14me12a}.
\bx

\subsect{9jl12e}{Proof of Theorem \ref{17me12c}}

We first show:

\prop{28me12a}{
Let $(X_1,\delta_1),(X_2,\delta_2),\ldots$ be complete metric spaces and let
$G$ be a group acting on each $X_k$, leaving
$\delta_k$ invariant ($k=1,2,\ldots$).
Then $(\prod_{k=1}^{\infty}X_k)/G$ is compact if and only if
$(\prod_{k=1}^tX_k)/G$ is compact for each $t$.
}

\pf
Necessity being direct, we show sufficiency.
We can assume that space $X_k$ has diameter at most $1/k$.
Let $A:=\prod_{k=1}^{\infty}X_k$, and let $d$ be the supremum metric
on $A$.
Then $d$ is $G$-invariant and
$\prod_{k=1}^{\infty}(X_k,\delta_k)$ is $G$-homeomorphic with $(A,d)$.
So it suffices to show that $(A,d)/G$ is compact.

As each $(X_k,\delta_k)$ is complete, $(A,d)$ is complete
(cf., e.g.,\ [3] Theorem XIV.2.5).
Hence, by Proposition \ref{2fe12a}, $(A,d/G)$ is complete.
So it suffices to prove that $(A,d/G)$ is totally bounded; that is, for each $\varepsilon>0$,
$A$ can be covered by finitely many
$d/G$-balls of radius $\varepsilon$.

Set $t:=\lfloor\varepsilon^{-1}\rfloor$.
Let $B:=\prod_{k=1}^tX_k$ and $C:=\prod_{k=t+1}^{\infty}X_k$,
with supremum metrics $d_B$ and $d_C$ respectively.
As $B/G$ is compact (by assumption), it can be covered
by finitely many $d_B/G$-balls of radius $\varepsilon$.
As $C$ has diameter at most $1/(t+1)\leq\varepsilon$, $A=B\times C$
can be covered by finitely many $d/G$-balls of radius $\varepsilon$.
\bx

\noindent
{\bf Proof of Theorem \ref{17me12c}:}
As each $(B_k,d_{R_k})$ is complete (Proposition \ref{28ja12e}),
by Proposition \ref{28me12a} it suffices to show that for each $t$,
$(\prod_{k=0}^t(B_k,d_{R_k}))/O(H)$ is compact.
Consider the Hilbert space $\prod_{k=0}^tH_k$,
and let $W:=\prod_{k=0}^tB_k$ and $R:=\prod_{k=0}^tR_k$.
Then the identity function on $W$ is
a homeomorphism from $(W,d_R)$ to $\prod_{k=0}^t(B_k,d_{R_k})$.
So it suffices to show that $(W,d_R)/O(H)$ is compact.
Now for each $n$, $R^n/O(H)$ is compact, as it is a continuous
image of $B(H)^m/O(H)$, with $m:=n(1+2+\cdots+t)$.
The latter space is compact, as it is a continuous image of the
compact space $B(\oR^m)^m$ (assuming $|C|=\infty$, otherwise $B(H)$
itself is compact).
So by Theorem \ref{13me12a}, $(W,d_R)/O(H)$ is compact.
\bx

\section*{References}\label{REF}
{\small
\begin{itemize}{}{
\setlength{\labelwidth}{4mm}
\setlength{\parsep}{0mm}
\setlength{\itemsep}{1mm}
\setlength{\leftmargin}{5mm}
\setlength{\labelsep}{1mm}
}
\item[\mbox{\rm[1]}] C. Borgs, J.T. Chayes, L. Lov\'asz, V.T. S\'os, K. Vesztergombi, 
Convergent sequences of dense graphs. I. Subgraph frequencies, metric properties
and testing,
{\em Advances in Mathematics} 219 (2008) 1801--1851.

\item[\mbox{\rm[2]}] J. Draisma, D. Gijswijt, L. Lov\'asz, G. Regts, A. Schrijver, 
Characterizing partition functions of the vertex model,
{\em Journal of Algebra} 350 (2012) 197--206.

\item[\mbox{\rm[3]}] J. Dugundji, 
{\em Topology},
Allyn and Bacon, Boston, 1966.

\item[\mbox{\rm[4]}] W. Fernandez de la Vega, R. Kannan, M. Karpinski, S. Vempala, 
Tensor decomposition and approximation schemes for constraint satisfaction problems,
in: {\em Proceedings of the 37th Annual {ACM} Symposium on Theory of Computing}
({STOC}'05),
pp. 747--754,
{ACM}, New York, 2005.

\item[\mbox{\rm[5]}] A. Frieze, R. Kannan, 
Quick approximation to matrices and applications,
{\em Combinatorica} 19 (1999) 175--220.

\item[\mbox{\rm[6]}] P.R. Halmos, 
{\em A Hilbert Space Problem Book --- Second Edition},
Springer, New York, 1982.

\item[\mbox{\rm[7]}] P. de la Harpe, V.F.R. Jones, 
Graph invariants related to statistical mechanical models:
examples and problems,
{\em Journal of Combinatorial Theory, Series B} 57 (1993) 207--227.

\item[\mbox{\rm[8]}] L. Lov\'asz, 
Graph homomorphisms: Open problems,
preprint, 2008.\\
\url{http://www.cs.elte.hu/~lovasz/problems.pdf}

\item[\mbox{\rm[9]}] L. Lov\'asz, 
{\em Large Networks and Graph Limits},
American Mathematical Society, Providence, R.I., 2012.

\item[\mbox{\rm[10]}] L. Lov\'asz, B. Szegedy, 
Limits of dense graph sequences,
{\em Journal of Combinatorial Theory, Series B} 96 (2006) 933--957.

\item[\mbox{\rm[11]}] L. Lov\'asz, B. Szegedy, 
Szemer\'edi's lemma for the analyst,
{\em Geometric and Functional Analysis} 17 (2007) 252--270.

\item[\mbox{\rm[12]}] T. Nishiura, 
Measure-preserving maps of $\oR^n$,
{\em Real Analysis Exchange} 24 (1998/9) 837--842.

\item[\mbox{\rm[13]}] A. Schrijver, 
Graph invariants in the edge model,
in: {\em Building Bridges --- Between Mathematics and Computer Science}
(M. Gr\"otschel, G.O.H. Katona, eds.), Springer, Berlin, 2008,
pp. 487--498.

\item[\mbox{\rm[14]}] B. Szegedy, 
Edge coloring models and reflection positivity,
{\em Journal of the American Mathematical Society}
20 (2007) 969--988.

\end{itemize}
}

\end{document}